\documentclass[11pt]{amsart}
\usepackage{times}
\usepackage[T1]{fontenc}
\usepackage{amscd,verbatim,amssymb}
\usepackage[all]{xy}
\usepackage{enumitem}
\usepackage{color}
\usepackage[colorlinks,linkcolor=blue,citecolor=blue,urlcolor=red]{hyperref}

\newcommand{\C}{\mathbf{C}}

\newcommand{\G}{\mathbb{G}}

\newcommand{\bH}{\mathbf{H}}

\renewcommand{\L}{\mathbb{L}}

\newcommand{\Q}{\mathbf{Q}}

\newcommand{\Z}{\mathbf{Z}}

\newcommand{\sA}{\mathcal{A}}
\newcommand{\sB}{\mathcal{B}}
\newcommand{\sC}{\mathcal{C}}
\newcommand{\sD}{\mathcal{D}}

\newcommand{\sG}{\mathcal{G}}

\newcommand{\sM}{{\operatorname{\mathbf{Mot}}}}

\newcommand{\fp}{\mathfrak{p}}
\newcommand{\fP}{\mathfrak{P}}

\newcommand{\un}{\mathbf{1}}

\newcommand{\dR}{{\operatorname{dR}}}
\newcommand{\cris}{{\operatorname{cris}}}

\newcommand{\et}{{\operatorname{\acute{e}t}}}

\newcommand{\GL}{{\operatorname{\mathbf{GL}}}}
\newcommand{\SL}{{\operatorname{\mathbf{SL}}}}

\newcommand{\Aut}{{\operatorname{Aut}}}
\newcommand{\Sm}{\operatorname{\mathbf Sm}}
\newcommand{\LMot}{\operatorname{\mathbf{LMot}}}
\newcommand{\Mot}{\operatorname{\mathbf{Mot}}}
\newcommand{\Ab}{\operatorname{\mathbf{Ab}}}

\newcommand{\Rep}{\operatorname{Rep}}

\newcommand{\Og}{{\operatorname{\mathbf{Og}}}}

\newcommand{\alg}{{\operatorname{alg}}}

\newcommand{\num}{{\operatorname{num}}}

\newcommand{\proj}{{\operatorname{proj}}}

\newcommand{\car}{\operatorname{char}}

\newcommand{\Spec}{\operatorname{Spec}}

\newcommand{\End}{\operatorname{End}}

\newcommand{\Ker}{\operatorname{Ker}}
\renewcommand{\Vec}{\operatorname{Vec}}

\renewcommand{\phi}{\varphi}
\renewcommand{\epsilon}{\varepsilon}

\newcommand{\inj}{\hookrightarrow}

\newcommand{\by}{\xrightarrow}

\newcommand{\iso}{\by{\sim}}

\renewcommand{\lim}{\varprojlim}

\newcounter{spec}
\newenvironment{thlist}{\begin{list}{\rm{(\roman{spec})}}%
{\usecounter{spec}\labelwidth=20pt\itemindent=0pt\labelsep=10pt}}%
{\end{list}}%

\newtheorem{Thm}{Theorem}
\newtheorem{Cor}[Thm]{Corollary}

\newtheorem{thm}{Theorem}[subsection]

\newtheorem{conjdrb}{De Rham-Betti conjecture}

\newtheorem{conjog}{Ogus conjecture}

\newtheorem{lemma}[thm]{Lemma}
\newtheorem{prop}[thm]{Proposition}
\newtheorem{cor}[thm]{Corollary}

\theoremstyle{definition}
\newtheorem{defn}[thm]{Definition}
\theoremstyle{remark}
\newtheorem{Rque}[Thm]{Remarks}
\newtheorem{rque}[thm]{Remark}

\renewcommand{\thesubsection}{\arabic{subsection}}

\numberwithin{equation}{subsection}

\begin{document}

\title[Fullness]{The fullness 
conjectures for products of elliptic curves}
\author{Bruno Kahn}
\address{CNRS, Sorbonne Université and Université Paris Cité, IMJ-PRG\\ Case 247\\4 place
Jussieu\\75252 Paris Cedex 05\\France}
\email{bruno.kahn@imj-prg.fr}
\address{Sorbonne Université, CNRS and Université Paris Cité, IMJ-PRG\\ Case 247\\4 place
Jussieu\\75252 Paris Cedex 05\\France}
\email{cyril.demarche@imj-prg.fr}
\date{October 7, 2024}
\maketitle

\hfill With an appendix by Cyril Demarche

\begin{abstract} We prove all conjectures of Yves Andr\'e's book \cite[Ch. 7]{andre} in the case of products of elliptic curves. The proofs given here are simpler and more uniform than the previous proofs in known cases.
\end{abstract}

%\enlargethispage*{40pt}

\subsection*{Introduction} Two of the most famous conjectures on a smooth projective variety $X$ over a suitable field $k$ are the Hodge and the Tate conjectures: $k=\C$ for the Hodge conjecture and $k$ finitely generated for the Tate conjecture. They are expounded for example by Yves André in his book on motives \cite[Ch. 7]{andre}, where he also describes two other similar conjectures, when $k$ is a number field:
\begin{itemize}
\item The ``de Rham-Betti'' conjecture  \cite[7.5.1.1]{andre}: it is related to a conjecture of Grothendieck on periods \cite[7.5.2.1]{andre} which can be traced back to \cite{groth}; see letter of Y. André to C. Bertolin in \cite{bertolin} for historical details.
\item The Ogus conjecture \cite[Introd. and \S 4]{ogus}; it is related to another conjecture of Grothendieck on algebraic solutions to differential equations \cite[p. 268]{ogus}.
\end{itemize}

For the readers' convenience, I recall these lesser-known conjectures in \S \ref{review}.

We shall be concerned here with the case where $X$ a product   of elliptic curves. If $k=\C$, the Hodge conjecture is known for $X$: this is attributed to Tate (unpublished) by Grothendieck in \cite[\S 3 c)]{grhodge}; a full proof was given by Imai in \cite{imai}.

Let $\ell$ be a prime number invertible in $k$, assumed to be finitely generated. The Tate conjecture for $\ell$-adic cohomology is known for $X$ in the following cases given in chronological order:

\begin{itemize}
\item $k$ is a number field (Imai \cite{imai});
\item $k$ is finite (Spie\ss\ \cite{spiess}).
\item $k$ is finitely generated over $\Q$: as pointed out in \cite{chen-vial}, this follows from \cite[Cor. 1.2]{lombardo}.
\end{itemize}

(In \cite{imai} and  \cite{lombardo}, Imai and Lombardo prove the Mumford-Tate conjecture for $X$, which implies the Tate conjecture from the Hodge conjecture.)

In each case, the result is stronger: the Hodge, or Tate, classes, are generated by those of degree $2$. The proofs, however, are different: for the Hodge and the Tate conjecture over a number field, Imai uses essentially a Tannakian argument plus results of Shimura-Taniyama while, over a finite field, Spie\ss' proof is obtained from an inequality on elliptic Weil numbers. On the other hand, the proof of the Hodge conjecture given by Imai involves a number of subcases and is especially delicate when dealing with elliptic curves with complex multiplication. 

This prompted me to look for a simpler and unified proof, which would also cover the two other conjectures. This was successful:

\begin{Thm}\label{t0} All the above-mentioned conjectures hold for $X$, in the strong sense that the cohomology classes coming from algebraic cycles are generated by those of degree $2$. In particular, the Tate conjecture holds for $X$ over any finitely generated field $k$. 
\end{Thm}

For the Ogus conjecture, the special case $k=\Q$, $X$ a power of a non CM elliptic curve is outlined in \cite[7.4.3.1]{andre}; it has inspired part of the proof here. Moreover, an assumption on $k$ in the first version of this article has now been dropped thanks to a theorem of Cyril Demarche (Theorem \ref{t2}).

The case of the de Rham-Betti conjecture in Theorem \ref{t0} implies part of a conjecture of Rohrlich: if all elliptic curves involved in $X$ have complex multiplication, then the multiplicative relations between its periods are generated by the ``obvious'' ones \cite[24.6.3.1]{andre}.  Another consequence of Theorem \ref{t0} is:

\begin{Cor}[see Proposition \ref{p3.1}]\label{c0} Let $X$ be a product of elliptic curves over any field $k$. Then the algebra of cycles modulo numerical equivalence on $X$ is generated in degree $1$.
\end{Cor}

\begin{Rque} a) Theorem \ref{t0} and its corollary extend to abelian varieties isogenous to products of elliptic curves (see  Proposition \ref{p3.1}).\\
b) To avoid a false impression, recall that there are known examples of abelian varieties $X$ such that the Hodge or the Tate conjecture hold for all powers of $X$ but where their ``strong form'' is false, e.g. \cite{shioda} or \cite[Ex. 1.8]{milne3}.\\
c) It it tempting to try and deduce the full Grothendieck period conjecture in the form of \cite[7.5.2.2]{andre} from Theorem \ref{t0}, say for a product $X$ of CM elliptic curves, reducing to the case of one such curve (Chudnovsky). However this is doomed to failure unless one knows something on the closure $Z$ of the canonical $\C$-point in the torsor $\fP$ of periods for $X$. Namely, a necessary condition is that $Z$ is a sub-torsor of $\fP$, meaning that its stabiliser in the corresponding Tannakian group has the same dimension as $Z$; conversely, this condition is inductively sufficient. Can one prove it?
\end{Rque}

\subsubsection*{Some words on the proof} The formalism developed here consists of two steps:

A) Assuming the ground field $k$ sufficiently large, it works  for any ``enriched realisation'' into a Tannakian category verifying certain axioms: see Theorem \ref{t1}. One then needs to check the axioms case by case: this is done in Section \ref{s5}. Here the key point is that each conjecture is known in codimension $1$ for abelian varieties: in the case of the Hodge conjecture, this is due to Lefschetz and Kodaira, for the Tate conjecture it is due to Tate, Zarhin and Faltings, and for the two other conjectures it was deduced by André in \cite{andre} from Wüstholz's analytic subgroup theorem for the de Rham-Betti conjecture and from results of Bost for the Ogus conjecture.

B) A descent argument.

To be able to use Tannakian arguments, one also needs the source category to be Tannakian. For this, and also to tackle the ``strong'' form of the conjectures conveniently, we use a category of ``Chow-Lefschetz motives'' (defined only for abelian varieties) introduced by  Milne in \cite{milne2} and developed in \cite{clmot}. 

As in \cite[7.6]{andre}, an important part of Step A) consists of a group-theoretic argument: it starts with the special case of one elliptic curve (and its powers), and uses a principle due to Goursat, Kolchin and Ribet to pass from there to the general case. The special case works well provided one knows that the Tannakian group attached to any elliptic curve is \emph{connected}. This is trivial for the Hodge conjecture and easy for the Tate conjecture; for the two other conjectures, I borrowed arguments from Yves Andr\'e (see \S\S \ref{s.period} and \ref{s.ogus}).

The general case is where elliptic curves with complex multiplication, whose motivic Galois groups are abelian, have rendered the Goursat-Kolchin-Ribet principle delicate in \cite{imai}. However, nobody seems to have used the full strengh of Kolchin's version of this principle \cite{kolchin}: his theorem is powerful enough to create a streamlined proof when the coefficients $K$ of the Weil cohomology involved are ``in good position with respect to the multiplications of the CM elliptic curves'': see Condition \ref{t1C} in Theorem \ref{t1}. This is true, in particular, when $K=\Q$, which is the case for the Hodge conjecture (where this approach trivialises Imai's arguments) and the de Rham-Betti conjecture.  For the two other conjectures, one needs more sophisticated arguments to get around this condition.

Step B) is much simpler than in the first version of this paper, relying on an elementary category-theoretic result (Proposition \ref{p6.1}).

As far as I have seen, the present method remains unfortunately very special to products of elliptic curves and not prone to generalisation. It raises nevertheless interesting questions about the generality of the \emph{result}. For example, let $X$ be an abelian variety of a type for which one of the conjectures is known ``in the strong sense'' (for $X$ and all its powers), e.g. one taken from the examples in \cite[A.7]{milne3} lifted to characteristic $0$. Can one prove the same for the other conjectures? At least, Corollary \ref{c3.1} shows that ``weak sense'' implies ``strong sense''  for these other conjectures.

This work was done in 2019, and was given a brief announcement in the algebraic geometry seminar of IMJ-PRG on June 20, 2019. %\footnote{Except that the Ogus conjecture was optimistically announced to be proven for $X$ over any number field, due to a incorrect weight argument.}. 
Since then, Kreutz, Shen and Vial have also proven the de Rham-Betti conjecture for products of elliptic curves in \cite{chen-vial}. Instead of Chow-Lefschetz motives, they use André's category of motivated motives and his ``Hodge = motivated'' theorem for abelian varieties to reduce to the Hodge conjecture.

\subsubsection{Notation}\label{not} We write $\Sm^\proj(k)=\Sm^\proj$ for the category of smooth projective varieties over a field $k$ and $\Ab(k)=\Ab$ for the category of abelian $k$-varieties and homomorphisms of abelian varieties.  We write $\Vec_K$ for the category of finite-dimensional vector spaces over a field $K$, and $\Rep_K(G)$ for the Tannakian category of finite-dimensional representations of an affine $K$-group $G$. If $\sC$ is an additive category, we write $\sC^{(\Z)}$ for the category of $\Z$-graded objects of $\sC$ with finite support. If $\sC$ is symmetric monoidal, we provide $\sC^{(\Z)}$ with the commutativity constraint given by the Koszul rule. A \emph{$\otimes$-functor} is a strong unital symmetric monoidal functor between unital symmetric monoidal categories.

%\enlargethispage*{20pt}

A full subcategory $\sD$ of an additive category $\sC$ is \emph{thick} if it is additive and stable under direct summands.

\subsection{Review of the de Rham-Betti and the Ogus conjectures}\label{review}

In both conjectures, $k$ is a number field; let $Y$ be a smooth projective $k$-variety. In the first case, we fix an embedding $k\inj \C$ and define a \emph{de Rham-Betti cycle of codimension $n$} as a pair $(\alpha,\beta)\in H^{2n}_\dR(Y/k)\times H^{2n}_B(Y,\Q)$ such that $\alpha\otimes 1\mapsto (2\pi i)^n\beta\otimes 1$ via the period isomorphism $H^{2n}_\dR(Y/k)\otimes_k \C\iso H^{2n}_B(Y,\Q)\otimes_\Q \C$. The cycle classes of any algebraic cycle of codimension $n$ yield a de Rham-Betti cycle and, conversely:

\begin{conjdrb}[\protect{\cite[7.5.1.1]{andre}}] Any de Rham-Betti cycle on $Y$ is algebraic.
 \end{conjdrb}

In the second case, we consider de Rham cohomology of $Y$ with extra structure: if $v$ is a finite unramified place of $k$ where $Y$ has good reduction, we have the Berthelot isomorphism \cite[3.4.2]{andre}
\[H^{2n}_\dR(Y_v/k_v)\iso H^{2n}_\cris(Y(v)/W(k(v)))\otimes_{W(k(v))} k_v\]
where $k_v$ (resp. $k(v)$) is the completion (resp. the residue field) of $k$ at $v$, $Y_v = Y\otimes_k k_v$ and $Y(v)$ is the special fibre of a smooth projective model of $Y$ at $v$. By transport of structure, the Frobenius automorphism of the right hand side provides the left hand side with an automorphism $\phi_v$ which is  semi-linear with respect to the absolute Frobenius of $k_v$. An \emph{Ogus cycle} is an element $\alpha$ of $H^{2n}_\dR(Y/k)$ such that, for almost all $v$, one has $\phi_v(\alpha) = q_v^n\alpha$, where $q_v=|k(v)|$. The cycle class of any algebraic cycle of codimension $n$ is an Ogus cycle and, conversely:

\begin{conjog}[\protect{\cite[7.4.1.2]{andre}}] Any Ogus cycle on $Y$ is algebraic.
\end{conjog}

\subsection{Background}\label{s2.1}
Let $k$ be a field, and let $\sM$ be the category of pure motives over $k$ modulo algebraic equivalence \cite[Ch. 4]{andre}, with coefficients in  a field $K$ of characteristic $0$. We shall use the notation $\sM(k)$ when it is necessary to specify $k$, but dispense from writing down coefficients $K$. We write $\L\in \sM$ for the Lefschetz motive. 

Let $\sB$ be a Tannakian category \cite[2.3]{andre} over $K$ (this means that $\End_\sB(\un)=K$). Surprisingly, I could not find a proof of the following result in the literature, while a corresponding result is available when $K$ is of characteristic $p>0$ \cite[Th. 6.1]{delignep}:

\begin{prop}\label{p2.2} In $\sB$, the tensor product and duals of two semi-simple objects $B,B'$ are semi-simple.
\end{prop}

\begin{proof} When $\sB$ is neutral, this follows from Chevalley's theorem \cite[p. 88]{chevalley}. In general, one may assume that $\sB$ is generated by $B$ and $B'$; by \cite[II, Rem. 3.10]{900}, there exists a fibre functor with values in a finite extension $L$ of $K$, that we may assume Galois of group $G$. The Tannakian category $\sB_{(L)}$ of \cite[pp. 155--156]{900} is neutralised by the canonical extension of this fibre functor along the inclusion $\sB\inj \sB_{(L)}$. Therefore it suffices to show that if an object $M\in \sB$ becomes semi-simple in $\sB_{(L)}$, it is semi-simple. Let $i:M'\inj M$ be a monomorphism. In $\sB_{(L)}$, the inclusion $i_{(L)}$ has a retraction $r$; the morphism $\frac{1}{|G|} \sum_{g\in G} gr g^{-1}$ is another, $G$-equivariant, retraction of $i_{(L)}$ which descends to a retraction of $i$.
\end{proof}

We can define a (generalised) Weil cohomology $\bH^*$ with values in $\sB$ just as in  \cite[3.3]{andre}: see \cite[VI.A.1]{saa} or \cite[4.2.1]{bvk}. It induces a $\otimes$-functor
\begin{equation}\label{d2.1}
\bH^*:\sM\to \sB^{(\Z)}
\end{equation}
(see \ref{not} for the notation). For $A\in \sB$ and $p\in \Z$, we write $A(p):=A\otimes \bH^2(\L)^{\otimes -p}$. 

If 
\[\omega:\sB\to \Vec_K\]
 is a neutral fibre functor, then $H^*=\omega^{(\Z)}\circ \bH^*$ is a Weil cohomology in the usual sense, and $\bH^*$ is an enrichment of $H^*$ in the sense of \cite[7.1.1]{andre}.

For $X\in \Sm^\proj$ and $r\ge 0$, we are interested in the condition
\begin{quote}
\begin{description}
\item[$F(X,r)$] the map 
\[A^r_\alg(X)\otimes K=\sM(\un,h(X)(r)) \to \sB(\un,\bH^{2r}(X)(r))\]
is surjective.
\end{description}
\end{quote}

The following is well-known, by a duality argument (cf. \cite[Ch. 7]{andre}):

\begin{lemma} \label{c1.1} Let $\sM[X]$ be the thick rigid $\otimes$-subcategory of $\sM$ generated by the motive $h(X)$ of some $X\in \Sm^\proj(k)$. Then Condition $F(X^n,r)$ for all $n,r\ge 0$ is equivalent to the fullness of $\bH^*$ restricted to $\sM[X]$.\qed
\end{lemma}

Suppose that $H$-homological equivalence agrees with numerical equivalence on $\sM[X]$. Then $H^*$ and therefore $\bH^*$ factor through the semi-simple $\otimes$-category $\sM_\num[X]$ \cite[4.3]{andre}. Moreover, the Künneth projectors of $X$ are algebraic (loc. cit., 5.4.2.1); in \eqref{d2.1}, after changing the commutativity constraints of $\sM_\num[X]$ as usual and those of $\sB^{(\Z)}$ and $\Vec_L^{(\Z)}$ by removing the signs of the Koszul rule, the composition of $\bH^*$ with the direct sum functor becomes symmetric monoidal and $H=\bigoplus H^*$ becomes a fibre functor. We may then study Condition $F(X,r)$ by Tannakian methods as in \cite[Ch. 7]{andre}, according to

\begin{prop}\label{p2.1} Let $\bH:\sA\to \sB$ be an exact $\otimes$-functor between Tannakian categories over $K$: in particular, $\bH$ is faithful by \cite[II, 1.19]{900}. Let $\omega:\sB\to \Vec_K$ be a neutral fibre functor and $H=\omega\circ \bH$. Suppose that $\sA$ has a $\otimes$-generator $M$ and that $\bH$ is essentially surjective. Then the Tannakian groups $G_H$ and $G_\omega$ of $H$ and $\omega$ are both subgroups of $\GL(H(M))$; in particular, $G_\omega\inj G_H$. Moreover, the following are equivalent:
\begin{thlist}
\item $G_\omega=G_H$;
\item $\bH$ is an equivalence of categories.
\end{thlist}
\end{prop}

\begin{proof}  The first claim is clear since a $\otimes$-automorphism of $H$ or $\omega$ is determined by its value on $M$ or $\bH(M)$. This said, (ii) $\Rightarrow$ (i) is trivial and (i) $\Rightarrow$ (ii) follows from Tannakian duality. \end{proof}

\begin{rque}\label{l2.1} If $\sA$ is semi-simple, the exactness of $\bH$ is automatic.
\end{rque}

Suppose that $k$ is of characteristic $0$, embeddable in $\C$, that $H^*$ is a classical Weil cohomology \cite[3.4]{andre}, and that $X$ is an abelian variety. Then homological equivalence agrees with numerical equivalence on $\sM[X]$ (loc. cit., 5.4.1.4) and we can use Proposition \ref{p2.1}. This approach has two drawbacks:

\begin{itemize}
\item it fails in positive characteristic;
\item in any characteristic, it does not account for the ``strong'' version of the fullness conjectures as explained in the introduction.
\end{itemize}

Our solution is to use Chow-Lefschetz motives, as explained next.

\subsection{Chow-Lefschetz motives} \label{s3}

When $k$ is algebraically closed, Milne defined in \cite{milne2} a rigid $\otimes$-category $\LMot$ of ``Lefschetz motives''\footnote{This terminology creates an ambiguity with the name ``Lefschetz motive'' given to $\L$; we hope this will not cause any confusion.} modelled on abelian varieties, for homological equivalence modulo  a Weil cohomology $H$ as above;   homological and numerical equivalences agree in this category, which therefore does not depend on the choice of $H$ and is abelian semi-simple. In \cite{clmot}, its definition was extended to other adequate equivalence relations and over not necessarily algebraically closed fields; as in \cite{milne2}, correspondences are given by sums of intersections of divisors. In this generality, numerical equivalence even agrees with \emph{algebraic} equivalence (loc. cit., Th. 3). In particular, if we work modulo algebraic equivalence then $\LMot$ maps to $\sM$, becomes Tannakian after changing the commutativity constraint as before, and $H$ induces a fibre functor on $\LMot$.  As above, we assume that the coefficients of $\Mot$ and $\LMot$ are $K$. We write $\iota:\LMot\to \Mot$ for the inclusion functor.  The additive functor $h^1:\Ab\otimes K\to \Mot$ of \cite[4.3.3]{andre} or \cite[\S 6.11]{zetaL} factors via $\iota$ through an additive functor $Lh^1(X):\Ab\otimes K\to \LMot$ \cite[Cor. 4.1 and Def. 4.4]{clmot}.

 \begin{lemma}\label{l3.1} Let $\LMot^1$ be the thick subcategory of $\LMot$ generated by the $Lh^1(X)$ where $X$ runs through $\Ab$, and define $\Mot^1$ similarly.  Then the functor $\iota^1: \LMot^1 \to \Mot^1$ induced by $\iota$  is an equivalence of categories.
 \end{lemma}
 
 \begin{proof}  First assume that $K=\Q$.  By  \cite[Th. 6.37]{zetaL}, $h^1$ is fully faithful, and essentially surjective by Poincaré's complete reducibility. Since  $\iota^1$ is  faithful, $Lh^1$ is then fully faithful. Therefore any idempotent endomorphism of some $Lh^1(X)$ comes from $\Ab\otimes \Q$, hence $Lh^1$ is an equivalence of categories and so is also $\iota^1$. This remains true after tensoring morphisms with $K$ and then taking pseudo-abelian envelopes, which concludes the proof since   $\LMot$ and $\Mot$ are pseudo-abelian.
  \end{proof}
 
If $X\in \Ab$, the category $\LMot[X]$  analogous to $\sM[X]$ (Lemma \ref{c1.1}) is $\otimes$-generated by the Chow-Lefschetz motive $Lh^1(X)$ \cite[Cor. 4.1]{clmot}, just as $\Mot[X]$ is $\otimes$-generated by $h^1(X)$ \cite[4.3.3]{andre}.  We have a string of $\otimes$-functors 
\begin{equation}\label{eq2.1}
\LMot[X]\by{\iota} \Mot[X] \by{\bH} \sB[X] \by{\omega} \Vec_K
\end{equation}
where $\sB[X]$ is the full Tannakian subcategory of $\sB$ \cite[2.3.5]{andre} generated by $\bH^1(X)$. We therefore have inclusions of Tannakian groups
\begin{equation}\label{eq2.2}
G_\omega(X)\inj G_{H\iota}(X)\inj \GL(H^1(X)). 
\end{equation}

\begin{prop}\label{p3.1} a) The categories $\LMot[X]$ and $\Mot[X]$  only depend on the isogeny class of $X$. The functor $\LMot^1[X]\to \Mot^1[X]$ is an equivalence of categories, where $\Mot^1[X]=\Mot^1\cap \Mot[X]$ and  $\LMot^1[X]=\LMot^1\cap \LMot[X]$. \\
b) The following are equivalent: referring to \eqref{eq2.1},
\begin{thlist}
\item  $\bH\iota$ is full;
\item the restriction of $\bH$ to $\Mot^1[X]$ is full and, for any $n\ge 0$, the graded algebra $\bigoplus_{r\ge 0} \sB(\un,\bH^{2r}(X^n)(r))$ is  generated in degree $1$. 
\item $F(X^n,r)$ holds for any $n,r\ge 0$ and the graded algebra $\bigoplus_{r\ge 0} A^r_H(X^n)$ is generated in degree $1$, where $A^*_H$ denotes algebraic cycles with $K$ coefficients modulo $H$-homological equivalence.
\end{thlist}
They imply
\begin{thlist}
\item[\rm(iii')] $F(X^n,r)$ holds for any $n,r\ge 0$ and the graded algebra $\bigoplus_{r\ge 0} A^r_\num(X^n)$ is generated in degree $1$, where $A^*_\num$ denotes algebraic cycles with $\Q$ coefficients modulo numerical equivalence.
\end{thlist}
c) The category $\sB[X]$ is semi-simple if and only if $\bH^1(X)$ is semi-simple. This happens if and only if $G_\omega(X)$ is reductive (not necessarily connected). \\
d) If $\bH^1(X)$ is semi-simple, Conditions (i) -- (iii) are equivalent to (iii') and imply that $\bH\iota$ is an equivalence of categories.\\
e) If $\bH^1(X)$ is semi-simple, Conditions (i) -- (iii) are also equivalent to
\begin{thlist}
\item[\rm(iv)] in \eqref{eq2.2}, $G_\omega(X)= G_{H\iota}(X)$.
\end{thlist}
\end{prop}

%\enlargethispage*{20pt}

\begin{proof} a) follows from Lemma \ref{l3.1} and its proof.  b) (i) $\Rightarrow$ (ii): the first part of (ii) follows from a), and the second part follows from Lemma \ref{c1.1} and the fact that the graded algebra $\bigoplus_{r\ge 0} \LMot(\un,Lh^{2r}(X^n)(r))$ is  generated in degree $1$ by definition of the morphisms in $\LMot$. 
(ii) $\Rightarrow$ (iii): the first claim is obvious and the second one then follows. (iii) $\Rightarrow$ (i):  again by Lemma \ref{l3.1}, the first part of (iii) imples that the restriction of $\bH\iota$ to $\LMot^1[X]$ is full, while its second part extends this fullness to all of $\LMot[X]$. (iii) $\Rightarrow$ (iii'): this is clear for algebraic cycles with $K$ coefficients modulo numerical equivalence, and we deduce the case with $\Q$ coefficients from \cite[3.2.7.1]{andre}. 

c) The first claim follows from Proposition \ref{p2.2}, and the second one from \cite[II, Rem. 2.28]{900}. 

d) (iii') $\Rightarrow$ (iii) follows from c) and \cite[7.1.1.1]{andre}. For the rest, in view of (i) it suffices to see that $\bH\iota$ is essentially surjective. But an object $H$ of $\sB[X]$ is a direct sum of direct summands of tensor constructions $\bH^1(X)^{\otimes p}\otimes \bH^1(X)^{*\otimes q}=\bH(Lh^1(X)^{\otimes p}\otimes Lh^1(X)^{*\otimes q})$, and we may assume that $H$ is a single such direct summand. By full faithfulness, the idempotent with image $H$ comes from an idempotent endomorphism of $Lh^1(X)^{\otimes p}\otimes Lh^1(X)^{*\otimes q}$. (This reasoning was also used in the proof of Lemma	 \ref{l3.1}; see \cite[Lemma 2.3.8]{bvk} for a general context.)

 e) In view of d), (i) $\iff$ (iv)  was seen in Proposition \ref{p2.1}. 
\end{proof}

\begin{cor}\label{c3.1} Let $\bH,\bH'$ be two enriched realisations. Suppose that $\bH^1(X)$ and ${\bH'}^1(X)$ are semi-simple and that  $F(X^n,r)$ holds for any $n,r\ge 0$ for both $\bH$ and $\bH'$. Then the conditions of Proposition \ref{p3.1} b) hold for $\bH$ if and only if they hold for $\bH'$.
\end{cor}

\begin{proof} Indeed, if (iii) holds for $\bH$, it implies (iii') for $\bH$ by Proposition \ref{p3.1} b), hence also for $\bH'$ by hypothesis, which in turn implies (iii) for $\bH'$ by hypothesis and Proposition \ref{p3.1} d).
\end{proof}

\begin{rque} One should beware that the last claim of (iii) in Proposition \ref{p3.1} b) becomes \emph{false} if one replaces homological equivalence with algebraic equivalence, already for $k=\C$ and $X$ the cube of the Fermat elliptic curve, cf. \cite[remarks after Th. 3]{clmot}.
\end{rque}

In all the cases given in the introduction, the first property of (ii) in Proposition \ref{p3.1} is verified (see Lemma \ref{l6.4}). We draw a consequence, weaker than (iv) in this proposition. Let $\bar K$ be an algebraic closure of $K$; if $G$ is a closed subgroup of $\GL_n$, write $\bar KG$ for the sub-$\bar K$-algebra of $M_n(\bar K)$ generated by $G(\bar K)$. Note that an irreducible representation of $G$ is also $G(\bar K)$-irreducible because $G(\bar K)$ is Zariski-dense in $G$; hence a semi-simple representation of $G$ is also $G(\bar K)$-semi-simple.

\begin{prop}\label{p3.3} Suppose that $\bH^1(X)$ is semi-simple and that the restriction of $\bH$ to $\Mot^1[X]$ is full. Then $\bar KG_\omega(X)=\bar K G_{H\iota}(X)$.
\end{prop}

\begin{proof} Both algebras are semi-simple since they admit the faithful semi-simple module  $\bar K\otimes_KH^1(X)$,  and they have the same centraliser by the fullness hypothesis and Lemma \ref{l3.1}. The conclusion now follows from the double centraliser theorem  \cite[\S 14, no 5, th. 5 a)]{bbki}.
\end{proof}

In the next section, we shall use

\begin{defn}\label{d4.1} Let $k_s$ be a separable closure of $k$, and $X_s=X\otimes_k k_s$. The field $k$ is \emph{large} relatively to $X$ if $Gal(k_s/k)$ acts trivially on $\End^0(X_s):=\End(X_s)\otimes \Q$, or equivalently if the injection $\End^0(X)\inj \End^0(X_s)$ is bijective.
\end{defn}

From now on we simplify the notation $G_{H\iota}(X)$ to $G_H(X)$. This group comes with a ``Tate'' character $t_X$, given by the $\otimes$-subcategory of $\LMot[X]$ generated by the Lefschetz motive $\L$\footnote{Note that $\L\in \LMot[X]$, as a direct summand of $Lh^2(X)$ given by the class of an ample line bundle.}. If $k$ is large relatively to $X$, as a special case of \cite[Th. 5]{clmot} we have an isomorphism
\begin{equation}\label{eq8}
G_H(X)\iso \prod_{\G_m} G_H(X_\alpha), 
\end{equation}
where the $X_\alpha$'s are the isotypic components of $X$ up to isogeny, and $\prod_{\G_m}$ denotes the fibre product over $\G_m$ with respect to the $t_{X_\alpha}$'s. (See \cite[(0.1)]{clmot} for the precise value of the $G_H(X_\alpha)$.) Note also that $k$ is large with respect to $X$ if and only if it is large with respect to $X_\alpha$ for each $\alpha$.

\begin{prop}\label{p3.2} In $\sB[X]$, $\bH^1(X)$ is semi-simple if and only if all the $\bH^1(X_\alpha)$ are semi-simple. If this holds and if $k$ is large with respect to $X$, the following are equivalent:
\begin{enumerate}[label={\rm (\arabic*)}]
\item Condition (iv) of Proposition \ref{p3.1} e) holds for $X$;
\item Condition (iv) of Proposition \ref{p3.1} e) holds for each $X_\alpha$ and the obvious homomorphism
\begin{equation}\label{eq8b}
G_\omega(X)\to \prod_{\G_m} G_\omega(X_\alpha),
\end{equation}
analogous to \eqref{eq8} is bijective.
\end{enumerate}
\end{prop}

\begin{proof} The first claim is obvious. (1) $\Rightarrow$ (2): the fullness of $\bH$ on $\LMot[X]$ implies its fullness on $\LMot[X_\alpha]$ for each $\alpha$, hence the first condition of (2) follows from Proposition \ref{p3.1} e); its second condition is now obvious. (2) $\Rightarrow$ (1) is also clear, thanks to \eqref{eq8}.
\end{proof}

\begin{rque}\label{r3.1} \eqref{eq8b} is a monomorphism, as follows from comparing it with \eqref{eq8}. Moreover, under the hypotheses of Proposition \ref{p3.2} the homomorphism $\pi_\alpha:G_\omega(X)\allowbreak\to G_\omega(X_\alpha)$ is faithfully flat by \cite[II, Rem. 2.29]{900}, since it corresponds to a full embedding of semi-simple Tannakian categories.
\end{rque}

\subsection{The case of a large ground field} We keep the notation of Sections \ref{s2.1} and \ref{s3}: in particular, $\sB$ is a Tannakian category over a field $K$ of characteristic $0$, $\omega$ is a neutral fibre functor on $\sB$ and $\bH=\bigoplus \bH^n$ where $\bH^*$ is a $\otimes$-functor as in \eqref{d2.1} , that we implicitly restrict to $\LMot$ (see \eqref{eq2.1}).

\begin{thm}\label{t1}  Let $X=\prod_{i\in I} E_i$ be a product of elliptic curves over $k$. Assume that $k$ is large relatively to $X$ and that
\begin{enumerate}[label={\rm (\Alph*)}]
\item \label{wh1}  $G_\omega(E_i)$ is connected for all $i$;
\item\label{t1A}  
$\bH^1(X)$ is semi-simple;
\item\label{t1B} 
the restriction of $\bH^1$ to $\Mot[X]$ is full (equivalently, the restriction of $\bH$ to $\Mot^1[X]$ is full);
\item\label{t1C} For all $(i,j)$ such that $E_i,E_j$ are non-isogenous and have complex multiplication (are ordinary if $\car k>0$), write $F_{ij}$ for the real quadratic subfield of  $\End(E_i)\otimes\End(E_j)\otimes\Q$. Then $F_{ij} \otimes_\Q K$ is a field.
\end{enumerate}
Then, in \eqref{eq2.1}, $\bH\iota:\LMot[X]\to \sB[X]$ is an equivalence of categories.
\end{thm}

Here are some comments on this statement. The assumption that $k$ be large is made to be able to use \eqref{eq8}. Note that the conclusion of Theorem \ref{t1} implies \ref{wh1}, \ref{t1A} and \ref{t1B}, because the groups $G_H(E_i)$ are connected  by \cite[Ex. 1]{clmot}  and $\LMot[X]$ is semi-simple. So the only ``mysterious'' condition is \ref{t1C}; its relevance will appear when we apply Kolchin's theorem below (Case (II) in 2)).

\begin{proof} Our strategy is to prove Condition (iv) of Proposition \ref{p3.1} by the method of Proposition \ref{p3.2}, i.e. to prove item (2) of this proposition. Thus there are two steps:

1) $X$ is a single elliptic curve. Note that the Tate character $t_X$ factors through $\det:\GL(H^1(X))\to \G_m$.  The group $G_H(X)$ is connected reductive by 
\cite[Ex. 1]{clmot}, and so is $G_\omega(E)$ by \ref{wh1} and \ref{t1A}. The connected reductive subgroups $G$ of $\GL_2$ which map surjectively to $\G_m$ by the determinant are $\GL_2$, $\G_m$ (its centre) and maximal tori. These subgroups are distinguished by the $\bar K$-algebras $\bar KG$,  and we conclude by Proposition \ref{p3.3}.

2) The general case. Write $M^1$ for $\Ker(M\to \G_m)$, for all groups $M$ appearing in the picture. For simplicity, write also $G=G_\omega(X)$ and $G_i=G_\omega(E_i)$. Let $J'=\{i\in I\mid G_i^1=\{1\}\}$: projecting onto $\prod_{i\in J-J'} G_i$, we may assume that $J'=\emptyset$.

It now suffices to prove that the inclusion $G^1\subseteq \prod_i G_i^1$ (see Remark \ref{r3.1}) is an equality.  Assume the contrary and write $\pi_i^1:G^1\to G_i^1$ for the projection induced by $\pi_i$: since $\pi_i$ is faithfully flat (Remark \ref{r3.1} again), so is $\pi^1_i$. By  \ref{wh1}, all $G_i$'s are connected, hence isomorphic to $\GL_2$ (no complex multiplication), $R_{\End(E_i)\otimes K/K} \G_m$ (complex multiplication/ordinary in positive characteristic) or $\G_m$ (supersingular). This  ensures that all $G_i^1$'s are quasi-simple ($K$-simple in the terminology of \cite{kolchin}), so we are in a position to apply Kolchin's result \cite[Theorem]{kolchin}. By this theorem  and the remarks following it, there are two possibilities:
\begin{enumerate}
\item[(I)] There exist two distinct indices $i\ne j$ such that $G_i^1, G_j^1$ are nonabelian (hence isomorphic to $\SL_2$), and an isomorphism $\phi:G_i^1\iso G_j^1$ such that the diagram
\[\xymatrix{
& G^1\ar[dl]_{\pi_i^1}\ar[dr]^{\pi_j^1}\\
G^1_i\ar[rr]^\phi&& G_j^1
}\]
commutes.
\item[(II)] There exist $l$ distinct indices $j(1), \dots ,j(l)$ with $l\ge 2$ and $G^1_{j(1)}$, \dots, $G^1_{j(l)}$, each commutative, and $l$ faithfully flat $K$-homo\-morph\-isms
$f_\lambda: G^1_{j(\lambda)}\to  G^1_{j(l)}$ ($1 \le \lambda\le  l$) such that $\prod_{1\le \lambda\le l} f_\lambda\circ \pi_{j(\lambda)}^1 = 1$.
\end{enumerate}

More precisely, (I) (resp. (II)) corresponds to (ii) (resp. (iii)) of loc. cit.\footnote{(i) is excluded by Step 1).}, except that (ii) is stated modulo finite subgroups of $G_i^1$ and $G_j^1$; but the remark in loc. cit., p. 1154 shows that one can get rid of these finite subgroups. In (II), the only thing we are going to use is the existence of a nontrivial homomorphism $f_\lambda:G^1_{j(\lambda)}\to  G^1_{j(l)}$ for some $\lambda< l$.

In Case (I), the image of $G^1$ in $G_i^1\times G_j^1$ must be the diagonal, hence the image of $G$ in $G_i\times_{\G_m} G_j$ must be the fibred diagonal $\Delta$. But, inside $M_4(\bar K)$, the two subalgebras $\bar K\Delta=\left\{\begin{pmatrix} u & 0\\0 &u\end{pmatrix}\mid u\in M_2(\bar K) \right\}$ and $\bar K(\GL_2\times_{\G_m}\GL_2)=\left\{\begin{pmatrix} u & 0\\0 &v\end{pmatrix}\mid u,v\in M_2(\bar K) \right\}$ are distinct, hence this case contradicts Proposition \ref{p3.3} in view of \eqref{eq8}.

In Case (II), recall that $1$-dimensional $K$-tori are classified by quadratic characters of $Gal(\bar K/K)$ via Cartier duality. Let $\chi_\lambda$ be the quadratic character corresponding to $G^1_{j(\lambda)}$: this is also the quadratic character corresponding to the étale $K$-algebra $\End(E_{j(\lambda)})\otimes K$. Suppose $\lambda<l$. Since $f_\lambda$ is not constant, we have $\chi_\lambda=\chi_l$ or equivalently $\chi_\lambda\chi_l=1$,  which contradicts \ref{t1C} since $\chi_\lambda\chi_l$ defines the étale $K$-algebra $F_{j(\lambda) j(l)}\otimes_\Q K$.
\end{proof}

\begin{rque}\label{r4.1} The proof of Theorem \ref{t1} shows that Condition \ref{t1C} is only needed when $X$ has at least two CM/ordinary factors.
\end{rque}

For the proof of the Tate conjecture (\S \ref{sTate}), we shall also need the following

\begin{prop} \label{p4.1} Assume that Conditions \ref{wh1}, \ref{t1A} and \ref{t1B} of Theorem \ref{t1} are satisfied. Write $X=X_1\times X_2$, where $X_1$ is the product of the CM factors of $X$ (ordinary factors in positive characteristic) and $X_2$ is the product of the other factors. Then the inclusion $G_\omega(X)\subseteq G_\omega(X_1)\times_{\G_m} G_\omega(X_2)$ is an equality.
\end{prop}

\begin{proof} Since the homomorphisms $G_\omega(X)\to G_\omega(X_i)$ are faithfully flat by Remark \ref{eq8b},  so are the homomorphisms $\pi_i:G^1_\omega(X)\to G^1_\omega(X_i)$. Let $N_i=\Ker \pi_i$: its image in $G^1_\omega(X_j)$ ($j\ne i$) is a normal subgroup. This time we apply Goursat's lemma \cite{goursat}: the image of $G^1_\omega(X)$ in $G^1_\omega(X_1)\times G^1_\omega(X_2)$ is the graph of an isomorphism $G^1_\omega(X_1)/N_2\iso G^1_\omega(X_2)/N_1$. But, by Step 1) in the proof of Theorem \ref{t1} and Remark \ref{r3.1}, $G^1_\omega(X_1)/N_2$ is a group of multiplicative type, and by Remark \ref{r4.1} $G^1_\omega(X_2)/N_1$ is a product of copies of $\SL_2$. Therefore these two groups are trivial, hence $G^1_\omega(X)$ contains $G^1_\omega(X_1)\times\{1\}$ and $\{1\}\times G^1_\omega(X_2)$ and the conclusion follows.
\end{proof}

For the proof of the Ogus conjecture, we shall need the following result of Cyril Demarche. 

\begin{thm}\label{t2} Let $F$ be a number field and let $\sG$ be a connected algebraic gerbe over $F$ \cite[II, Appendix]{900}\footnote{Here ``connected algebraic'' means that, for any extension $E$ of $F$ such that $\sG(E)\neq \emptyset$ and for any $X\in \sG(E)$, the affine $E$-group $G_X$ such that $G_X(K)=\Aut(X_K)$ for any extension $K/E$ is of finite type and connected. This is equivalent to the condition ``with connected linear band'' in Proposition \ref{pdemarche}.}. Then, for any finite (field) extension $R/F$, there exists a finite extension $K/F$, linearly disjoint from $R/F$, such that $\sG(K)\neq \emptyset$.
\end{thm}

\begin{proof} See Proposition \ref{pdemarche} in the appendix. 
\end{proof}

\subsection{Descent} Suppose that we have an enrichment as in Definition \ref{d2.1}, but don't assume $k$ large. We would like to reduce Theorem \ref{t0} to Theorem \ref{t1}; it is the aim of this section.

%\enlargethispage*{30pt}

\subsubsection*{The idea}   Let $\Gamma$ be a finite group and let $f^*:\sA\to \sA'$ be a $\Gamma$-equivariant functor, where $\sA$ and $\sA'$ are categories provided with a (pseudo-)action of $\Gamma$, this action being trivial on $\sA$. Then $f^*$ induces a functor from $\sA$ to the category $\sA'[\Gamma]$ of descent data with respect to the action of $\Gamma$; we say that $f^*$ \emph{has descent} if this functor is an equivalence of categories \cite[end of \S 1.1]{stacks}.  Let $l/k$ be a finite Galois extension and $\Gamma=Gal(l/k)$. As observed in the proof of  \cite[Prop. 7.1]{stacks},  it follows from the definition of $\LMot$ \cite[Def. 4.4]{clmot} and the proof of  Theorem 5 in loc. cit., §5.5 that $f^*:\LMot(k)\to \LMot(l)$ has descent, where $f:\Spec l\to \Spec k$ is the structural morphism. This fact restricts to the categories $\LMot(k)[X]$ and $\LMot(l)[X_l]$ for any $X\in \Ab(k)$.

Suppose now that we have a naturally commutative diagram of rigid abelian $\otimes$-categories and $\otimes$-functors
\begin{equation}\label{eq4.30}
\begin{gathered}
\xymatrix{
\LMot(l) \ar[r]^{\bH(l)}&\sB(l)\\
\LMot(k)\ar[u]_{f^*} \ar[r]^{\bH}&\sB\ar[u]_{f_\sB^*}
}
\end{gathered}
\end{equation}
in which $\bH(l)$ restricts to an equivalence of categories on the rigid subcategories generated by $Lh^1(X)$ and its images, as in the previous sections. If $f_\sB^*$ also has descent and $\bH(l)$ is $\Gamma$-equivariant, we get the same conclusion for $\bH$. 

The problem in this strategy is that, in practice, the category $\sB(l)$ which is used to formulate a given fullness conjecture fails to carry a $\Gamma$-action, so that the above argument does not make sense. As an example, for the Tate conjecture we use $\sB=\Rep_{\Q_\ell}(G_k)$ (continuous representations), where $G_k=Gal(k_s/k)$ is the absolute Galois group of $k$ relative to some separable closure $k_s$. But the action of $G_k$ by conjugation on $\Rep_{\Q_\ell}(G_l)$, where $l$ is supposed to be a subextension of $k_s/k$, does not factor through its quotient $\Gamma$. The case of de Rham-Betti realisations is similar.

Fortunately there is a catch-all solution to this issue, which is given by the following proposition.

\begin{prop}\label{p6.1} Consider a naturally commutative diagram of categories and functors
\begin{equation}\label{eq4.32}
\begin{gathered}
\xymatrix{
\sA' \ar[r]^{\bH'}&\sB'\ar[r]^{\rho'}&\sC'\\
\sA\ar[u]_{f^*} \ar[r]^{\bH}&\sB\ar[u]_{f_\sB^*}\ar[r]^\rho& \sC\ar[u]_{f_\sC^*}
}
\end{gathered}
\end{equation}
in which a finite group $\Gamma$ acts on $\sA'$ and $\sC'$. We assume that 
\begin{itemize}
\item $\rho'\bH'$ is $\Gamma$-equivariant, as well as $f^*$ and $f_\sC^*$ for the trivial action of $\Gamma$ on $\sA$ and $\sC$;
\item $\rho'$ and $f_\sB^*$ are faithful; 
\item the functor $\sA\to \sA'[\Gamma]$ induced by $f^*$ is fully faithful (in particular, $f^*$ is faithful).
\end{itemize}
 If $\bH'$ is fully faithful, so is $\bH$.
\end{prop}

\begin{proof} This is a standard categorical diagram chase, that we make explicit for the benefit of the reader. For simplicity, we reason as if \eqref{eq4.32} were strictly commutative.

Let $A,A'\in \sA$ and let $w,w'\in \sA(A,A')$ be such that $\bH(w)=\bH(w')$. Then $\bH'f^*(w)=\bH'f^*(w')$, hence $w=w'$ since $f^*$ and $\bH'$ are faithful. This proves that $\bH$ is faithful.

Let now $u\in \sB(\bH(A),\bH(A'))$. By hypothesis, there exists $v\in \sA'(f^*A,f^*A')$ such that $\bH'(v)=f_\sB^* (u)$. I claim that $gv=v$ for any $g\in \Gamma$; by the faithfulness of $\rho'$ and $\bH'$, it suffices to show this after applying $\rho'\bH'$. But $\rho'\bH'(v)= \rho'f_\sB^* (u)=f_\sC^*\rho(u)$, hence 
\[\rho'\bH'(gv)=g\rho'\bH'(v)=gf_\sC^*\rho(u)=f_\sC^*\rho(u)=\rho'\bH'(v)\] 
where we used the $\Gamma$-equivariance of $\rho'\bH'$. Using now the hypothesis on $f^*$, we find that $v=f^*w$ for some $w\in \sA(A,B)$. Hence 
\[f_\sB^*\bH(w)=\bH'f^*(w)=\bH'(v)=f_\sB^*(u)\]
and $\bH(w)=u$ by the faithfulness of $f_\sB^*$. This proves that $\bH$ is full.
\end{proof}

\begin{rque}\label{r5.1}
Proposition \ref{p6.1} does not give the essential surjectivity of $\bH$ even if we assume it for $\bH'$.  Theorem \ref{t0} does not state such essential surjectivity. Yet, one can deduce it from Proposition \ref{p3.1} d) by the semi-simplicity of $\bH^1(X)$ (Lemma \ref{l6.4}) if one wishes.
\end{rque}

\subsection{Proof of Theorem \ref{t0}}\label{s5} First we have:

\begin{lemma}[{\cite[th. 7.1.7.5]{andre}}]\label{l6.4} Conditions \ref{t1A} and \ref{t1B} of Theorem \ref{t1} hold for each conjecture of \cite[Ch. 7]{andre}.\qed
\end{lemma}

(This was the starting point of this paper.)

We now go case by case by using Proposition \ref{p6.1} to reduce to the case where $k$ is large enough, so that we may apply Theorem \ref{t1}.

\subsubsection{The Hodge conjecture \cite[7.2]{andre}}\label{sHodge} Here $k=\C$, $H$ is Betti cohomology, $\sB$ is the category $\mathbf{Hdg}$ of pure polarisable $\Q$-Hodge structures. Since $\C$ is algebraically closed, there is no descent issue. Condition \ref{wh1} holds because Mumford-Tate groups are connected, cf. \cite[7.1.2.1 1)]{andre}, and Condition \ref{t1C} is trivial because $K=\Q$.

\subsubsection{The Tate conjecture \cite[7.3]{andre}} \label{sTate} Here $k$ is finitely generated over its prime field, $H$ is $\ell$-adic cohomology for some prime number $\ell$ invertible in $k$, $\sB$ is the category of continuous representations of $G_k=Gal(k_s/k)$ on finite-dimensional $\Q_\ell$-vector spaces, and $K=\Q_\ell$.  

We may replace the latter by the equivalent Tannakian category $R_\ell(k)$ of continuous representations of the absolute Galois groupoid $\Pi_k$, as in \cite[2.2.3]{adjoints}. The corresponding functor $\bH$ sends $h(X)$, for $X\in \Sm^\proj(k)$, to the functor from $\Pi_k$ to $\Vec_{\Q_\ell}$ which sends a separable algebraic closure $k_s$ of $k$ to $H^*_\et(X_{k_s},\Q_\ell)$. If $l/k$ is a finite Galois extension with group $\Gamma$, then $\Gamma$ acts on $R_\ell(l)$ by the functoriality of the étale fundamental groupoid. We may then apply Proposition \ref{p6.1} with $\sA=\LMot(k)$, $\sA'=\LMot(l)$, $\sB=\sC=R_\ell(k)$ and $\sB'=\sC'=R_\ell(l)$. The faithfulness of $f_\sB^*$ is trivial and the functor $\bH'$ corresponding to $\bH$ ``over $l$'' is $\Gamma$-equivariant. To prove Theorem \ref{t0} in the case of the Tate conjecture, we may therefore assume $k$ large relatively to $X$, and further enlarge it if needed. 

For $M\in \sB$, $G_\omega(M)$ is the Zariski closure of the action of $G_k$ on $\omega(M)$ \cite[7.1.3]{andre}; it follows that the composition $G_k\to G_\omega(M)(k)\to \pi_0(G_\omega(M))$ is surjective, which implies Condition \ref{wh1} up to enlarging $k$ to a finite Galois extension $k_\ell$, possibly depending on $\ell$. 

By quadratic reciprocity, Condition \ref{t1C} holds for a set of prime numbers in a suitable union of arithmetic progressions depending on $X$, which completes the proof of the Tate conjecture (in its strong form) for these $\ell$.  We are now going to prove it for \emph{any} $\ell$ by avoiding Condition \ref{t1C}, in several steps. 
\bigskip

\noindent\emph{6.2.1. $k$ is finite.} By the above, there exists at least one prime $\ell_0$ such that the conclusion of Theorem \ref{t1} holds for $\ell_0$-adic  cohomology. In particular, the $\ell_0$-adic Tate conjecture holds for $X$, and $H^1_{\ell_0}(X)$ is semi-simple. By Condition (c) of \cite[Th. 2.9]{tate}, the same then holds for any $\ell$. Therefore the conclusion of Theorem \ref{t1} holds for any $\ell$ thanks to Corollary \ref{c3.1}. 
\bigskip

\noindent\emph{6.2.2. $k$ is a number field.} Write $X=X_1\times X_2$, where $X_1$ is the product of the CM factors of $X$ and $X_2$ is the product of the other factors. By Theorem \ref{t1}, we know that $G_\omega(X_2)= G_H(X_2)$ (because Condition  \ref{t1C} is empty for $X_2$), and  by Proposition \ref{p4.1} we know that  $G_\omega(X)\iso G_\omega(X_1)\times_{\G_m} G_\omega(X_2)$. Therefore, by \eqref{eq8}, it suffices to prove that $G_\omega(X_1)= G_H(X_1)$. Choose an embedding $\sigma$ of $k_s$ into $\C$, and let $H_B:\LMot\to \Vec_\Q$ be the corresponding Betti realisation, which factors as $\LMot\by{\bH_B} \mathbf{Hdg}\by{\omega'}\Vec_\Q$ where $ \mathbf{Hdg}$ is as in \S\ref{sHodge}. Recall M. Artin's comparison isomorphism $H\simeq H_B\otimes_\Q \Q_\ell$. We have $G_{\omega'}(X_1)=G_{H_B}(X_1)$ by \S\ref{sHodge}  and $G_{\omega'}(X_1)\otimes_\Q\Q_\ell=G_\omega(X_1)$ by \cite[p. 371]{imai} or  \cite[\S 4]{chiafu}, which concludes the proof.
\bigskip

\noindent\emph{6.2.3. The general case.} Let $k_0$ be the field of constants of $k$. Write $X=X_1\times X_2$ as above, where ``CM'' means ``ordinary'' in positive characteristic.  As in 6.2.2., we reduce to the case $X=X_1$. Note that $X$ is then defined over a finite extension of $k_0$, hence taking $k$ sufficiently large we may assume that $X$ is actually defined over $k_0$. Then the isomorphism $G_\omega(X)\iso G_H(X)$ over $k_0$ remains valid over $k$, because the left and the right hand sides do not change under such extension of scalars:  the surjectivity of $G_k\to G_{k_0}$ implies $\LMot(k_0)[X]\iso \LMot(k)[X]$ on the right hand side, while on the left hand side it implies that $G_k$ and $G_{k_0}$ have same Zariski closure in $\GL(H^1(X))$. 

\subsubsection{The de Rham-Betti conjecture \cite[7.5]{andre}}\label{s.period} Here $k$ is a number field embedded in $\C$, $H$ is Betti cohomology, $\sB$ is the category $\Vec_{k,\Q}$ of \cite[7.1.6]{andre}. Here again, let $l/k$ be a finite Galois extension with group $\Gamma$. We apply Proposition \ref{p6.1} with $\sA=\LMot(k)$, $\sA'=\LMot(l)$, $\sB=\Vec_{k,\Q}$, $\sC'=\Vec_k$, $\rho$ the forgetful functor,  and similarly for $\sB',\sC'$ and $\rho'$ replacing $k$ with $l$. As in \S \ref{sTate}, the faithfulness of $f_\sB^*$ is trivial and so is that of $\rho'$; $\rho'\bH'$ is the de Rham realisation which is $\Gamma$-equivariant. We may therefore assume $k$ to be large.

Since $K=\Q$, Condition \ref{t1C} is trivial. It remains to prove Condition \ref{wh1} up to extending $k$. Let $\bar \Q\subset \C$ be the set of algebraic numbers, so that $k\subset \bar \Q$. We start with a lemma:

\begin{lemma} Write $G_N$ for the Tannakian group of the full Tannakian subcategory $\langle N\rangle\subset \Vec_{k,\Q}$ generated by an object  $N=(W,V,\varpi)$. For an extension $l$ of $k$ contained in $\bar \Q$, write $N_l=(l\otimes_k,W,V,\varpi)$ for the image of $N$ under the obvious functor $\Vec_{k,\Q}\to \Vec_{l,\Q}$. Then there exists a finite extension $l/k$ such that the natural homomorphism $G_{N_{\bar \Q}}\to G_{N_l}$ is an isomorphism.
\end{lemma}

\begin{proof} The affine group schemes $G_N$ and $G_{N_l}$ are algebraic by \cite[II, Prop. 2.20 (b)]{900}, and the homomorphism $G_{N_l}\to G_N$ is a closed immersion because the criterion of \cite[II, Prop. 2.21 (b)]{900} is trivially verified by the description of the objects of $\langle N\rangle$ and $\langle N_l\rangle$ as subquotients of direct sums of tensor constructions \cite[I, 3.1a]{900}. Therefore the inverse system $(G_{N_l})$, where $l$ runs through the finite subextensions of $\bar \Q/k$, is stationary, and it suffices to show that its inverse limit is $G_{N_{\bar \Q}}$. It suffices to verify this on $R$-points for any $\Q$-algebra $R$ (or even for $R=\bar \Q$). But this is clear, since $\langle N_{\bar \Q}\rangle$ is the $2$-colimit of the $\langle N_l\rangle$'s.
\end{proof}

We are  now left to show that any finite quotient $G$ of the Tannakian group of $\Vec_{\bar \Q,\Q}$ is reduced to $\{1\}$. 

The following simple argument was kindly communicated by Y. Andr\'e. Let $(W,V, i)$ be the object of $\Vec_{\bar \Q,\Q}$ corresponding to a representation of $G$. Then $i$ is defined over $\bar \Q$ (indeed, $W, V$  and their tensor constructions define a  torsor under $G$, and $i$ defines a complex point, hence a $\bar \Q$-point  of this finite torsor). The choice of a basis of $V$ then identifies $(W,V, i)$ with a sum of copies of the unit object in $\Vec_{\bar \Q,\Q}$.

\subsubsection{The Ogus conjecture}\label{s.ogus}  Here $k$ is a number field embedded in $\C$, $H=H_\dR$ is de Rham cohomology, $\sB$ is the category $\Og(k)$ of \cite[7.1.5]{andre} and $K=\Q$ \cite[p. 72, footnote (4)]{andre}. Since $\Og(l)$ carries a natural action of $\Gamma$ and the Ogus realisation is $\Gamma$-equivariant, we proceed directly with this category as in \S\ref{sTate}, and may assume $k$ large.

The fibre functor  $\omega:\Og(k)\to \Vec_k$ is not neutral (unless $k=\Q$), but it becomes so after extending it to the Tannakian category over $k$ $\Og(k)_{(k)}$  of \cite[pp. 155-156]{900} (notation: $\omega_k$). 

We prove Condition \ref{wh1} for the corresponding Tannakian group $G_{\omega_k}(E)$ of an elliptic factor $E$ of $X$, by proving directly the equality $G_{\omega_k}(E)=G_H(E)$.
This follows from the argument in \cite[7.4.3.1]{andre}. A detailed sketch: choose a prime $\fp$ of $k$, unramified over $\Q$ and with good and ordinary reduction for $E$, so that Tsuji's $p$-adic period isomorphism  $H^1_\dR(E/k)\otimes_k B_{pst}\simeq H^1_\dR(E_{k_\fp}/k_\fp)\otimes_{k_\fp} B_{pst}\simeq H^1_\et(\bar E_{k_\fp},\Q_p)\otimes_{\Q_p} B_{pst}$ \cite[3.4.2]{andre} is defined. The action of $G_{\omega_k}(E)$ on the left group is compatible with the action of $Gal(\bar k_\fp/k_\fp)$ on the right group. If $E$ is not CM, Corollary 1 of \cite[p; IV-44]{ladic} says that the Lie algebra of  the latter action is a Borel subalgebra of $\mathfrak{gl}_2(\Q_p)$. Hence $G_{\omega_k}(E)=\GL_2(=G_H(E))$, since it is reductive. If $E$ has complex multiplication, Corollary 2 of loc. cit. says that this Lie algebra is a split Cartan algebra; hence $G_{\omega_k}(E)$ contains a maximal torus; but then it must be equal to the maximal torus $G_H(E)$ \cite[Ex. 1]{clmot}.

Let now $\sG$ be the gerbe associated to $\sB[X]$. If $R$ denotes the compositum of the fields $F_{ij}$ appearing in Condition \ref{t1C}, there exists by Theorem \ref{t2} another fibre functor $\omega':\Og(k)\to \Vec_K$ where $K$ is a number field linearly disjoint from $R$. Replacing $\Og(k)$ by $\Og(k)_{(K)}$, we get a neutral fibre functor $\omega'_K$; since gerbes are locally connected, Condition \ref{wh1} for $\omega_k$ implies Condition \ref{wh1} for $\omega'_K$. By construction, Condition \ref{t1C} is satisfied, hence $\bH$ becomes fully faithful after extending scalars from $\Q$ to $K$, and therefore it is fully faithful.

\subsection*{Acknowledgements} I thank Yves Andr\'e for several discussions around this paper, Giuseppe Ancona for pointing out \cite{chen-vial}, Cyril Demarche for kindly writing up the appendix and the referee for several comments which helped improve the exposition.

\renewcommand{\thesubsection}{\Alph{subsection}}

\setcounter{subsection}{0}

\subsection{Appendix: sections of algebraic gerbes}\label{sapp}\

\hfill Cyril Demarche

\begin{lemma} \label{lem appendix}
Let $F$ be a number field, $E, F_1, \dots, F_r$ be finite field extensions of $F$. For all $1 \leq i \leq r$, let $\alpha_i \in \textup{Br}(F_i)$.

Then there exists a finite field extension $K/F$, linearly disjoint from $E$ and from the $F_i$'s over $F$, and such that $\alpha_i$ vanishes in $\textup{Br}(K \otimes_F F_i)$ for all $i$.
\end{lemma}

\begin{proof}
For all $i$, let $S_i$ denote the finite set of places $v$ of $F_i$ such that $\alpha_{i,v} \neq 0$ in $\textup{Br}(F_{i,v})$. For all $i$ and all $v \in S_i$, denote by $n_{i,v}$ the order of $\alpha_{i,v}$ in $\textup{Br}(F_{i,v})$.

Let $S$ be the (finite) set of places $v$ of $F$ for which there exists an $i$ and a place $w \in S_i$ dividing $v$. For each $v \in S$, let $n_v$ denote the lcm of the $n_{i,w} [F_{i,w}:F_v]$ for all $i$ and $w \in S_i$ such that $w$ divides $v$.

Let $E' / F$ denote the Galois closure of the compositum of the fields $E$, $F_1$, \dots, $F_r$, and denote its Galois group by $G$. By Chebotarev theorem, for all $g \in G$, there exists a place $v_g$ of $F$ outside $S$ and not dividing $2$ such that $E'/F$ is unramified at $v_g$ and the Frobenius at $v_g$ lies in the conjugacy class of $g$. One can assume that $v_g \neq v_h$ if $g \neq h$.

Define $S'$ to be the union of $S$ and all the places $v_g$ for $g \in G$.

By the Grunwald-Wang theorem, there exists a (cyclic) field extension $K/F$ such that for all $v \in S$, all places $w$ of $K$ above $v$, the local degree $[K_w:F_v]$ is divisible by $n_v$, and for all $v \in S' \setminus S$, the extension $K/F$ splits completely at $v$.

Then by construction the extensions $K$ and $E'$ are linearly disjoint over $F$, and for all $i$, the image of $\alpha_i$ in $\textup{Br}(K \otimes_F F_i)$ vanishes at all places of $K \otimes_F F_i$ by a restriction-corestriction argument, hence it vanishes globally.
\end{proof}

We refer to \cite{Bor} and \cite{DLA} for the notions of gerbes, non-abelian $H^2$ and $k$-kernels.

\begin{prop}\label{pdemarche}
Let $k$ be a number field and $\mathcal{G}$ be an étale (algebraic) gerbe over $k$, with connected linear band (or $k$-kernel) $L$. Let $E$ be a finite extension of $k$.

Then there exists a finite extension $K/k$ linearly disjoint from $E/k$ such that $\mathcal{G}(K) \neq \emptyset$.
\end{prop}

\begin{proof}
Consider the class $\alpha$ of $\mathcal{G}$ in the non-abelian cohomology set $H^2(k,L)$ and its image $\alpha_F \in H^2(F,L)$. Let $T$ be the $k$-torus associated to the $k$-band $L$ (see \cite[1.7 and 6.1]{Bor}  for instance) and $\alpha' \in H^2(k,T)$ be the image of $\alpha$.

By Proposition 6.5 in \cite{Bor}, for any totally imaginary finite field extension $K/k$, $\mathcal{G}(K) \neq \emptyset$ if and only if $\alpha'_K = 0$ in $H^2(K,T)$.

So we are reduced to find a totally imaginary finite extension $K/k$ linearly disjoint from $E$ over $k$, such that $\alpha'$ vanishes in $H^2(K,T)$.

First, there exists a totally imaginary finite field extensions $F/k$ that is linearly disjoint from $E$ over $k$ (one can work as in Lemma \ref{lem appendix}).

So we are reduced to find a finite extension $K/F$ linearly disjoint from $E' := EF$ over $F$, such that $\alpha'$ vanishes in $H^2(K,T)$.

There exists an exact sequence of $k$-tori (a flasque resolution of $T$ for instance, see \cite[Proposition 1.3]{CTS})
\[0 \to S \to P \to T \to 0\]
such that $P$ is quasi-trivial, i.e. isomorphic to a finite product of tori of the shape $R_{k_i/k}(\mathbf{G_m})$ for some finite field extensions $k_i / k$. 
The field $F$ is totally imaginary, hence $H^3(F,S) = 0$ by \cite[Chapter I, Corollary 4.21]{ADT}. Therefore, the map $H^2(F,P) \to H^2(F,T)$ is surjective. Let $\beta' \in H^2(F,P)$ be a lift of the class $\alpha' \in H^2(F,T)$.

Writing the quasi-trivial $F$-torus $P_F$ as $\prod_{i=1}^r R_{F_i/F}(\mathbf{G_m})$ for some finite field extensions $F_i/F$, we get that 
\[H^2(F,P) \cong \bigoplus_{i=1}^r \textup{Br}(F_i) \, .\]

To conclude, we apply the lemma to the elements $(\beta_i)_{1 \leq i \leq r} \in \bigoplus_{i=1}^r \textup{Br}(F_i)$ corresponding to the class $\beta' \in H^2(F,P)$ and we get a finite extension $K/F$, linearly disjoint of $E'$ and the $F_i$'s over $F$, such that the image of $\beta'$ in $H^2(K,P) \cong \bigoplus_{i=1}^r \textup{Br}(K \otimes_F F_i)$ vanishes. Then $K$ and $E$ are linearly disjoint over $k$ and the proof is complete.
\end{proof}


\begin{thebibliography}{12}
\bibitem{andre} Y. Andr\'e Une introduction aux motifs (motifs purs, motifs mixtes,
p\'eriodes), Panoramas et synth\`eses {\bf 17}, SMF, 2004.
\bibitem{bvk} L. Barbieri-Viale, B. Kahn {\it Universal Weil cohomology}, preprint, 2024, \url{https://arxiv.org/abs/2401.14127}.
\bibitem{bertolin} C. Bertolin {\it Third kind elliptic integrals and $1$-motives} (with a letter of Y. André and an appendix by M. Waldschmidt), J. pure appl. Algebra {\bf 224} (2020), 106396.
\bibitem{Bor} M. Borovoi {\it Abelianization of the second nonabelian Galois cohomology}, Duke Math. J. {\bf 72} no.1, 217–239 (1993).
\bibitem{bbki} N. Bourbaki Algèbre, ch. 8: {\it Modules et anneaux semi-simples}, Springer, 2012.
\bibitem{chevalley} C. Chevalley Théorie des groupes de Lie, tome III, Hermann 1954.
%\bibitem{clifford} A. H. Clifford,  {\it Representations induced in an invariant subgroup},  Ann. of Math. {\bf 38} (1937), 533--550.
\bibitem{CTS} J.-L. Colliot-Thélène and J.-J. Sansuc {\it Principal homogeneous spaces under flasque tori; applications.} J. Algebra {\bf 106}, No. 1, 148-205 (1987). 
\bibitem{900} %P. Deligne and J. S. Milne {\it Tannakian categories}, {\it in} 
P. Deligne, J. S. Milne, A. Ogus and K. Shih Hodge cycles, motives, and Shimura varieties, Lect. Notes in Math. {\bf 900}, Springer-Verlag, 1982, 101--228. 
\bibitem{delignep} P. Deligne {\it Semi-simplicité de produits tensoriels en caractéristique p}, Invent. Math. {\bf 197} (2014), no. 3, 587--611.
\bibitem{DLA} C. Demarche and G. Lucchini Arteche {\it Le principe de Hasse pour les espaces homogènes : réduction au cas des stabilisateurs finis}, Compositio Math. {\bf 158} (8), 1568-1593 (2019).
%\bibitem{den-mur} C. Deninger and J. P. Murre {\it Motivic decomposition of abelian schemes and the Fourier transform},  J. Reine Angew. Math. {\bf 422} (1991), 201--219.
\bibitem{goursat} E. Goursat {\it Sur les substitutions orthogonales et les divisions régulières de l'espace}, Ann. Sci. ÉNS {\bf 6} 1889), 9--102
\bibitem{groth} A. Grothendieck {\it On the de Rham cohomology of algebraic varieties}, Publ. Math. IHÉS {\bf 29} (1966) 96--103.
\bibitem{grhodge} A. Grothendieck {\it Hodge's general conjecture is false for trivial reasons}, Topology {\bf 8} (1969), 299--303.
\bibitem{imai} H. Imai {\it On the Hodge group of some abelian varieties}, Kodai Math. Sem.
Rep. {\bf 27} (1976), 367--372.
%\bibitem{cell} B. Kahn {\it Équivalences rationnelle et numérique sur certaines variétés de type abélien sur un corps fini}, Ann. Sci. \'Ec. Norm. Sup. {\bf 36} (2003), 977--1002.
%\bibitem{ss} B. Kahn {\it On the semi-simplicity of Galois actions}, Rendiconti Sem. Mat. Univ. Padova {\bf 112} (2004), 97--102.
\bibitem{adjoints} B. Kahn {\it Motifs et adjoints}, Rend. Sem. mat. Univ. Padova {\bf 139} (2018), 77--128.
\bibitem{zetaL} B. Kahn Zeta and $L$-functions of varieties and motives, LMS Lect. Notes Series {\bf 462}, Cambridge Univ. Press, 2020.
\bibitem{clmot} B. Kahn {\it Chow-Lefschetz motives},  Indag. Math. 2024 (Murre memorial volume), \url{https://doi.org/10.1016/j.indag.2024.04.007}.
\bibitem{stacks} B. Kahn {\it Galois descent for motivic theories}, preprint, 2024, \url{https://arxiv.org/abs/2312.01825}.
\bibitem{kolchin} E. R. Kolchin {\it  Algebraic groups and algebraic dependence}, Amer. J. Math. {\bf 90} (1968),  1151--1164.
\bibitem{chen-vial} T. Kreutz, M. Shen,  C. Vial {\it de Rham--Betti classes on products of elliptic curves over a number field are algebraic}, preprint, 2023, \url{https://arxiv.org/abs/2206.08618}.
%\bibitem{lieb} D. Lieberman {\it Numerical and homological equivalence of algebraic cycles on Hodge manifolds}, Amer. J. Math. {\bf 90} (1968), 366--374.
\bibitem{lombardo} D. Lombardo {\it On the $\ell$-adic Galois representations attached to nonsimple abelian varieties}, Ann. Inst. Fourier (Grenoble) {\bf 66} (2016), 1217--1245. 
%\bibitem{mcl} S. Mac Lane Categories for the working mathematician, (2nd. ed.), Springer, 1998.
%\bibitem{milne} J.S. Milne {\it Lefschetz classes on abelian varieties}, Duke Math. J. {\bf 96} (1999), 639--675.
\bibitem{milne2} J.S. Milne {\it Lefschetz motives and the Tate conjecture}, Compositio Math.
{\bf 117} (1999), 47--81.
\bibitem{milne3} J.S. Milne {\it The Tate conjecture for certain abelian varieties over finite fields},  
Acta Arith. {\bf 100} (2001), 135--166.
\bibitem{ADT} J.S. Milne {\it Arithmetic duality theorems}, 2nd ed. Charleston, SC: BookSurge, LLC. viii, 339 p. (2006).  
%\bibitem{mov} {\it An "explicit" description of cocomma-categories ?}, Mathoverflow discussion, \break \url{https://mathoverflow.net/questions/247280/an-explicit-description-of-cocomma-categories/247311#247311}.
\bibitem{ogus} A. Ogus {\it Hodge cycles and crystalline cohomology}, {\it in} Hodge cycles, motives and Shimura varieties, Lect. Notes in Math. {\bf 900}, Springer, 1982, 357--414.
%\bibitem{pohlmann} H. Pohlmann {\it Algebraic cycles on Abelian varieties of complex multiplication type}, Ann. of Math. {\bf 88} (1968), 161--180.
\bibitem{saa} N. Saavedra Rivano Catégories tannakiennes, Lect. Notes in Math. {\bf 265}, Springer, 1972.
\bibitem{shioda} T. Shioda {\it Algebraic cycles on abelian varieties of Fermat type}, Math. Ann. {\bf 258} (1981), 65--80.
\bibitem{spiess} M. Spie\ss\ {\it Proof of the Tate conjecture for products of elliptic curves over finite fields}, Math. Ann. {\bf 314} (1999), 285--290. 
\bibitem{ladic} J.-P. Serre Abelian $l$-adic representations and elliptic curves, Benjamin, 1968.
\bibitem{tate} J. Tate {\it Conjectures on algebraic cycles in $l$-adic cohomology}, {\it in} Motives (Seattle, WA, 1991), Proc. Sympos. Pure Math., {\bf 55} (1), Amer. Math. Soc., 1994, 71--83.
\bibitem{chiafu} C.-F. Yu {\it A note on the Mumford-Tate conjecture for CM abelian varieties}, Taiwanese J. Math.
{\bf 19} (2015), 1073--1084.
\end{thebibliography}
\end{document}